\documentclass{ifacconf}

\usepackage{soul} 

\usepackage{amsmath,amssymb,mathtools}
\usepackage{graphicx}      
\usepackage{natbib}        
\usepackage{xcolor}
\usepackage{enumitem}
\usepackage[integrals]{wasysym}
\usepackage{booktabs}
\usepackage{array}

\makeatletter
\long\def\@makefigurecaption#1#2{%
  \vskip 5\p@\small
  \par\hangindent2em\@fignumfont{#1}. #2\par\vskip5pt
}
\long\def\@maketablecaption#1#2{%
  \vskip 5\p@\small
  \par\hangindent2em\@fignumfont{#1}. #2\par\vskip5pt
}
\makeatother

\renewcommand{\Re}{\mathbb{R}}

\newcommand{\NNb}{\mathbb{N}}

\newcommand{\nn}{{n\times n}}
\newcommand{\bd}{\mathrm{bd}}
\newcommand{\nint}{\mathrm{int}}
\newcommand{\diff}{\mathrm{d}}
\newcommand{\realp}{\mathrm{Re}}

\newcommand{\conv}{\mathrm{conv}}

\newcommand{\calE}{\mathcal{E}}
\newcommand{\calF}{\mathcal{F}}

\newcommand{\calL}{\mathcal{L}}

\newcommand{\calR}{\mathcal{R}}

\newcommand{\Sub}{\mathrm{Sub}}

\newcommand{\extended}[1]{{#1}}
\newcommand{\short}[1]{}

\newenvironment{talign}
    {\let\displaystyle\textstyle\align}
    {\endalign}
\newenvironment{talign*}
    {\let\displaystyle\textstyle\csname align*\endcsname}
    {\endalign}

\begin{document}

\setstcolor{red}  
\sethlcolor{brown} 

\begin{frontmatter}

\title{Bounds on set exit times of affine systems, using Linear Matrix Inequalities\thanksref{footnoteinfo}}

\thanks[footnoteinfo]{The work of the first author is supported by a FRIA (F.R.S.--FNRS) fellowship.}

\author[First]{Guillaume O.\ Berger}
\author[Second]{Maben Rabi}

\address[First]{Universit\'e catholique de Louvain (UCLouvain), B-1348 Louvain-la-Neuve, Belgium (e-mail: guillaume.berger@uclouvain.be)}
\address[Second]{\O stfold University College, NO-1757 Halden, Norway (e-mail: maben.rabi@hiof.no)}

\begin{abstract}
Efficient computation of trajectories of switched affine systems becomes possible, if for any such hybrid system, we can manage to efficiently compute the sequence of switching times.
Once the switching times have been computed, we can easily compute the trajectories between two successive switches as the solution of an affine ODE.
Each switching time can be seen as a positive real root of an analytic function, thereby allowing for efficient computation by using  root finding algorithms. These algorithms require a finite interval, within which to search for the switching time.
In this paper, we study the problem of computing upper bounds on such switching times, and we restrict our attention to stable time-invariant affine systems.
We provide semi-definite programming models to compute upper bounds on the time taken by the trajectories of an affine ODE to exit a set described as the intersection of a few generalized ellipsoids.
Through numerical experiments, we show that the resulting bounds are tighter than bounds reported before, while requiring only a modest increase in computation time.
\end{abstract}

\begin{keyword}
Hybrid and switched systems modelling, reachability analysis, verification and abstraction of hybrid systems
\end{keyword}

\end{frontmatter}

\section{Introduction}\label{sec-introduction}

Simulation of hybrid systems has been quite an important topic of research in recent years; see, e.g., \citet{carloni2006languages,goebel2009hybrid,schweiger2019anempirical} for surveys.
Hybrid systems are generally characterized by the presence of piecewise continuous dynamics.
A major challenge for the simulation of such systems is thus to determine the switching times of trajectories, that is, when they hit a \emph{switching surface} (boundary between two continuous components of the vector flow) \citep{cremona2019hybrid}.
Most hybrid system simulation tools rely on step refinement algorithms for reliably detecting the event of the trajectory crossing a boundary \citep{esposito2001accurate,wang2015astate,copp2016azerocrossing,farkas2019adaptive}.
These methods are variations of classical ODE solvers, which also use step refinement techniques to compute accurate solutions of continuous ODE flows.
But these methods handle both linear flows and nonlinear flows in the same way; nor do they discriminate between linear, ellipsoidal, or general nonlinear and non-convex switching surfaces.

On the other hand, some classes of systems allow for a \emph{piecewise analytic} expression of their solutions, making the use of numerical-integration-based ODE solvers superfluous.
This is the case for instance for  autonomous \emph{switched affine systems,} for which the trajectories of each individual mode can be expressed using elementary transcendental functions of the type: $e^{at}$, $\sin(\omega t)$, $\cos(\omega t)$ and $t^k$.
These systems constitute a paradigmatic class of hybrid and cyber-physical systems, and appear naturally in many engineering applications or as abstractions of more complicated systems; see, e.g., \citet{heemels2001equivalence,bell2010thecontinuous,legat2020stability}.
The recent work of \citet{rabi2020piecewise} proposes to use the piecewise analytic expression of the trajectories and to combine it with root-finding algorithms for the detection of boundary crossing.
This allows to benefit from the long-standing maturity of root-finding solvers for transcendental equations \citep{boyd2014solving}, in order to compute fast and accurate solutions for these systems.

The main challenge of the above approach is to produce reliable time intervals for the crossing of a switching surface by the trajectory, as root-finding algorithms generally require a bounded interval in which to search for potential roots of the analytic function.
In \citet{rabi2020piecewise}, upper bounds on the crossing time are computed using Lyapunov functions.
The idea is that if the current mode of the system is stable and its equilibrium lies away from the switching surface, then after some time the trajectory will enter an invariant region not intersecting the switching surface.
Quadratic Lyapunov functions are used to compute such invariant regions and upper bounds on the time for reaching them.
The approach in \citet{rabi2020piecewise} relies on arbitrarily restricting the Lyapunov equation to have a simple right-hand side, as in $A^\top\! P+PA=-I$.
A geometric interpretation of the resulting Lyapunov function leads to upper bounds on the switching time that grow quadratically with the distance of the initial point to the boundary.

In this work, we push further the above approach by combining it with convex optimization techniques to reduce its conservativeness.
The idea is to use semi-definite programs (SDPs) to compute Lyapunov-like functions, and invariant regions that are not necessarily centered at the equilibrium.
The SDPs shall be formulated to provide tight upper bounds on the crossing times.
The upper bounds obtained in this way are optimal within the considered framework of Lyapunov-like functions, and grow either quadratically or logarithmically with the distance of the initial point to the boundary, depending on the SDP formulation.
When our method is used to compute switched system trajectories, each continuous mode is associated with a corresponding SDP, whose solution can be computed off-line.

We study two such optimization models for the computation of upper bounds on the switching time of switched affine systems, and compare them in terms of computation time and tightness of the bounds.
We also show improvements over the bounds of \citet{rabi2020piecewise}.
Our comparisons are empirical and are based on a large set of random matrices.
We observe that more complex models indeed lead to substantially tighter bounds while not requiring a significant increase of the computation time (which can be performed off-line if needed).
We focus on the comparison of the different models for computing upper bounds on the switching times, and we leave for further work the implementation and tests of these methods in a standalone software for the simulation of switched affine systems.

The paper is organized as follows.
In Section~\ref{sec-problem}, we introduce the problem of interest.
In Section~\ref{sec-upper-bounds}, we present the different models for the computation of upper bounds on the switching time.
In Section~\ref{sec-comparisons}, we present the numerical experiments for the comparisons of the different models.

{\bfseries Notation.}
For $N\in\NNb_{>0}$, we let $[N]=\{1,\ldots,N\}$.
For $A\in\Re^\nn$, $\sigma(A)$ denotes the \emph{stability margin} of $A$, that is, $\sigma(A)=\min\,\{-\realp(\lambda) : \text{$\lambda$ is an eigenvalue of $A$}\}$.

\section{Problem formulation}\label{sec-problem}

Consider an affine system $\dot{x}(t)=Ax(t)+b$, with $A\in\Re^\nn$ strictly stable and $b\in\Re^n$.
Let $\calR\subseteq\Re^n$ be a closed convex region with nonempty interior.
For $x_0\in\calR$, we define the \emph{escaping time from $\calR$} of the trajectory starting at $x_0,$ by
\[
t_*(x_0,\calR;A,b)=\min\,\{t\geq0:\xi(t,x_0)\in\bd\,\calR\},
\]
where $\xi(t,x_0)$ is the trajectory of the affine system given by $A$ and $b$, starting from $x_0$.
If $\xi(t,x_0)\in\nint\,\calR$ for all $t\geq0$, we let $t_*(x_0,\calR;A,b)=0$. 

If the system is a switched affine system, $\dot{x}=A_ix+b_i$ when $x\in\calR_i$, then the escaping times of $\dot{x}=A_ix+b_i$ from the regions $\calR_i$ gives the switching times of the trajectory.
The convention that $t_*(x_0,\calR_i;A_i,b_i)=0$ if the trajectory remains in $\calR_i$ is motivated by the problem of computing the switching time using root-finding algorithms, for which we seek upper bounds as small as possible on the interval in which the first crossing, \emph{if it exists}, occurs.

To simplify the notation, we will assume in the rest of the paper that each mode is \emph{linear}, i.e., has the form $\dot{x}=A_ix$.
Thus we will focus on the problem of finding upper bounds on the escaping time of the LTI system $\dot{x}=Ax$, where $A$ is strictly stable.
This conversion of stable affine systems into stable linear systems can be done, \emph{without loss of generality}, by translating the state $x$ and the region $\calR$ around the equilibrium point $\bar{x}=-A^{-1}b$.

We also make the following assumptions on the region $\calR$:
\begin{itemize}[itemsep=5pt]
\item $\calR$ is described as the intersection of a finite set of  ellipsoids%
\footnote{In this paper, the term \emph{ellipsoid} is understood in the broad sense, meaning that an ellipsoid can be degenerate (or unbounded in some directions).
However, an ellipsoid is always assumed to be closed and to have nonempty interior.}\!%
:
\[\textstyle
\calR=\bigcap_{i=1}^N \calE_i,
\]
where each $\calE_i$ is an ellipsoid.
Note that this includes the case of a polyhedral region $\calR$.
\item $\calR$ is included in a set described as the convex hull of a finite set of ellipsoids or singletons:
\[\textstyle
\calR\subseteq\conv\,\{\calF_1,\ldots,\calF_M\},
\]
where each $\calF_i$ is an ellipsoid or a singleton.
\item The origin does not belong to the boundary of $\calR$.
\end{itemize}

\begin{rem}
In our analysis, we assume that the origin is not on the boundary of $\calR$ and that $A$ is strictly stable.
This implies that for every $x_0\in\calR$ there is a finite time $T_*\geq0$ such that for all times $t\geq T_*$, the state $\xi(t,x_0)$ is bounded away from the boundary of $\calR$.
On the other hand, if the origin was on the boundary of $\calR$, or if $A$ was unstable or marginally stable, then there could be situations in which the trajectories of the system tend to the boundary of $\calR$ but never reach it, or tend to the boundary of $\calR$ and cross it an infinite number of times (e.g., convergence with rotation); so that this more general class of matrices would be very difficult to handle for the problem of escaping time estimation, at least with the techniques described in this work.
Some techniques, based on the eigenvalue decomposition, were proposed in \citet{rocca2015exploiting} to exploit the rotation of the trajectories to detect reachability of LTI systems.
However, since even with these techniques it is not guaranteed to obtain a finite upper bound on the escaping time, we preferred to focus here on stable affine systems.
\end{rem}

\begin{rem}
We assume that the convex set $\calR$ can be enclosed in the convex hull of a set of ellipsoids or singletons.
This includes the two extreme cases: when $\calR$ is enclosed in a single ellipsoid and when it is enclosed in the convex hull of a finite set of points.
When the number of points increases, the second case can provide arbitrarily accurate enclosures of $\calR$.
However, for some regions (e.g., the hypercube), the number of points needed to describe the convex hull grows exponentially with the dimension.
In these cases, it may be beneficial to consider enclosures described by a few ellipsoids (instead of many points), even if it provides less accurate enclosures of $\calR$.
\end{rem}

\section{Upper bounds on switching times}\label{sec-upper-bounds}

In this section, we present two frameworks, inspired by Lyapunov theory, to compute upper bounds on the escaping time of stable LTI systems.
By computing Lyapunov-like functions that have guaranteed rates of decrease inside $\calR$, we can find invariant regions not intersecting the boundary of $\calR$ and obtain upper bounds on the time for the trajectories to enter the invariant region.
To compute these functions, we follow the approach of quadratic Lyapunov theory, and formulate their existence as the feasibility of a convex optimization program in which the constraints on the functions are encoded as \emph{Linear Matrix Inequalities} (LMIs).

\subsection{Quadratic functions as decision variables}\label{ssec-quad-functions}

A \emph{quadratic function} (on $\Re^n$) is a function of the form $V(x)=x^\top Qx + 2b^\top x+c$ with $Q=Q^\top\in\Re^\nn$, $b\in\Re^n$ and $c\in\Re$.
A quadratic function $V(x)=x^\top Qx + 2b^\top x+c$ is \emph{convex} if $Q\succeq0$.
Given a quadratic function $V$ and a scalar $s\in\Re$, we let $\Sub_s(V)$ be the \emph{sublevel set} $\{x\in\Re^n:V(x)\leq s\}$.
Note that any ellipsoid (possibly degenerate) can be described as $\Sub_0(V)$ for some convex quadratic function $V$.
A quadratic function $V$ is \emph{nonnegative} if $V(x)\geq0$ for all $x\in\Re^n$.
This will be denoted by $V\geqq0$.
We also use the notation $V_1\geqq V_2$ for $V_1-V_2\geqq0$, where $V_1,V_2$ are two quadratic functions.

In our framework, we will use quadratic functions as \emph{decision variables}.
More precisely, we will define optimization problems where the objective is to find unknown quadratic functions that have to satisfy several constraints.
Because quadratic functions are parameterized by $Q$, $b$ and $c$, finding them will in fact amount to compute the associated $Q$, $b$ and $c$.
Moreover, the constraints on the quadratic functions will always have the form $V_1\geqq V_2$.
The proposition below states that this type of constraints can be encoded as LMIs, so that the unknown quadratic functions can be computed using Semi-Definite Programming.

\begin{prop}\label{pro-quad-form-positive}
The quadratic function $V(x)=x^\top Qx + 2b^\top x+c$ is nonnegative ($V\geqq0$) if and only if $\left[\begin{array}{cc} Q & b \\ b^\top & c \end{array}\right]\succeq0$.
\end{prop}

\begin{pf}
This follows from $V(x)=[x\,;1]^\top \left[\begin{array}{cc} Q & b \\ b^\top & c \end{array}\right] [x\,;1]$.
See for instance \citet{boyd2004convex}.\hfill\qed
\end{pf}

Given a matrix $A$ and a quadratic function $V$, the \emph{Lie derivative} of $V$ along the vector field $x\mapsto Ax$ is defined as $\calL_AV(x) = V'(x)Ax$.
Observe that if $V(x)=x^\top Qx + 2b^\top x+c$, then $\calL_AV(x) = x^\top(A^\top Q+QA)x+2b^\top Ax$; in particular, $\calL_A V$ depends linearly on $V$.

In the next subsections, we introduce optimization models to compute upper bounds on the set escaping times of stable linear systems.
For the sake of readability, these models are presented using quadratic functions as variables.
From the above discussion, these models can be reformulated as Semi-Definite Programs in the associated variables $Q$, $b$ and $c$, thereby allowing for efficient computation of the solution, using for instance interior-point algorithms \citep{bental2001lectures,boyd2004convex}.

In the following, we let $F_1,\ldots,F_M$ and $E_1,\ldots,E_N$ be convex quadratic functions whose $0$-sublevel sets are equal to the ellipsoids describing $\calR$ (see Section~\ref{sec-problem}); i.e., $\Sub_0(F_i)=\calF_i$ for all $i\in[M]$, $\Sub_0(E_i)=\calE_i$ for all $i\in[N]$.

\subsection{Upper bound on escaping time. Case I: $0\in\nint\,\calR$}\label{ssec-bound-escaping-time-inside}

This is the case where the origin is in the interior of the region $\calR$.
Consider the following optimization program.

\emph{Decision variables:} convex quadratic functions $V$ and $W$, and scalars $\lambda_1,\ldots,\lambda_M\geq0$, $\mu_1,\ldots,\mu_N\geq0$, $\nu_1,\ldots,\nu_N\geq0$ and $r\in\Re$;
\emph{Objective and constraints:}
\begin{subequations}\label{eq-EscapeIn}
\begin{talign}
&\!\!\makebox[10mm][l]{max}  r \label{eq-EscapeIn-obj}\\
&\!\!\makebox[10mm][l]{s.t.} V \leqq 1 + \lambda_i F_i, \quad \forall\,i\in[M], \label{eq-EscapeIn-include-R0}\\
&\!\!\makebox[10mm][l]{}     V \geqq r + \mu_i E_i, \quad \forall\,i\in[N], \label{eq-EscapeIn-inside-R1}\\
&\!\!\makebox[10mm][l]{}     \calL_AW \leqq 0, \label{eq-EscapeIn-W-deriv}\\
&\!\!\makebox[10mm][l]{}     W \geqq \nu_i E_i, \quad \forall\,i\in[N], \label{eq-EscapeIn-W-inside-R1}\\
&\!\!\makebox[10mm][l]{}     \calL_AV \leqq G(V) - W, \label{eq-EscapeIn-deriv}
\end{talign}
\end{subequations}
where $G(V)$ is equal to either $-1$ (leading to quadratically increasing bounds), or to $-2\gamma V$ for some fixed parameter $\gamma>0$ (leading to logarithmically increasing bounds).

\begin{figure}
\centering
\includegraphics[width=0.99\linewidth]{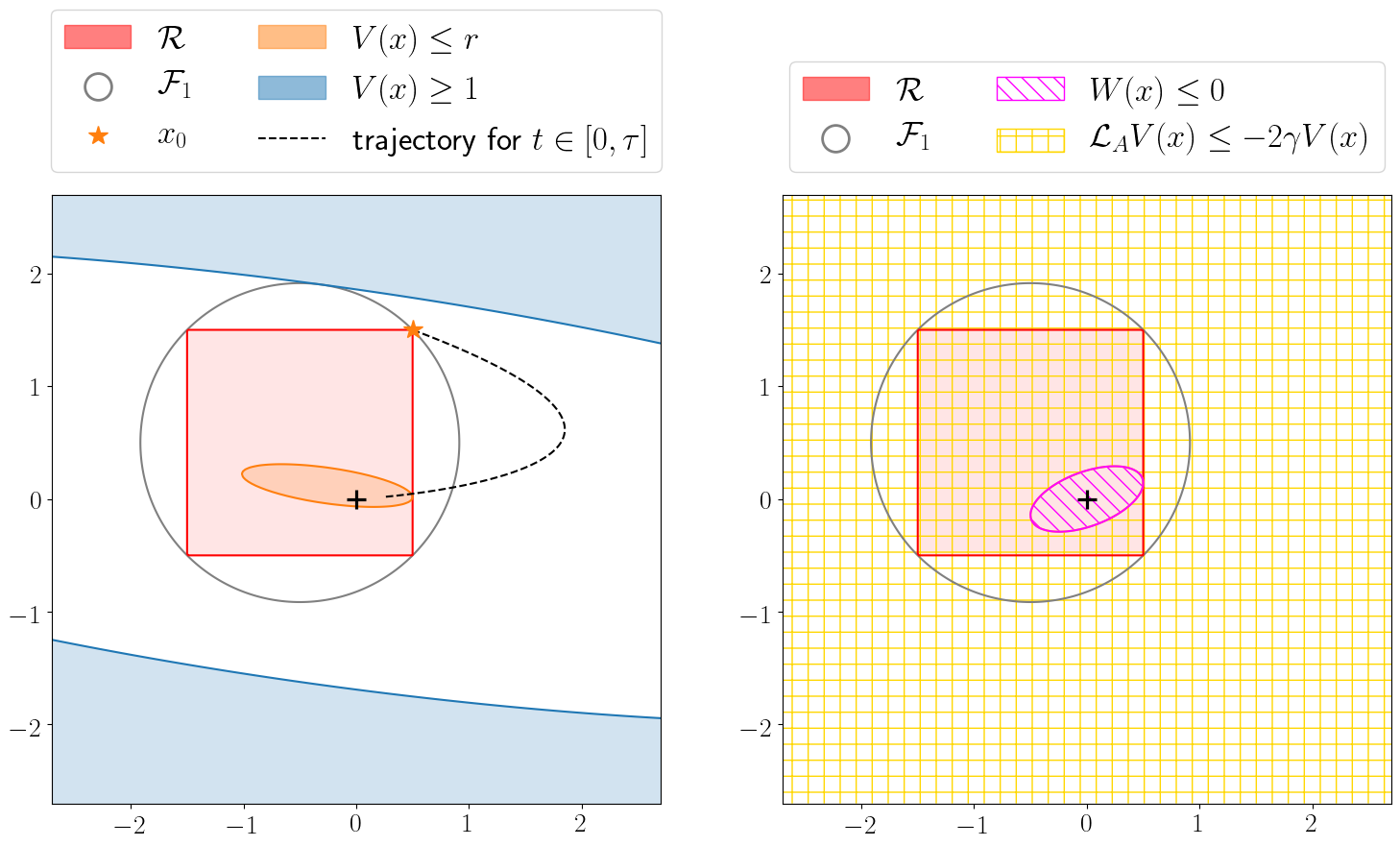}
\caption{Illustration of Example~\ref{exa-EscapeIn-simple}.
We verify that $\calR\subseteq\Sub_1(V)$, $\Sub_r(V)\subseteq\calR$, and $\calL_A V(x)\leq G(V(x))$ whenever $W(x)\geq0$.}
\label{fig-EscapeIn-simple}
\end{figure}

{\bfseries Explanations.} The goal of these constraints is to find an ellipsoidal invariant region, defined by $\Sub_0(W)$ where $W$ is a convex quadratic function, that is included in $\calR$.
This is enforced by the constraints \eqref{eq-EscapeIn-W-deriv} (implying that the Lie derivative of $W$ is nonpositive everywhere, so that any sublevel set of $W$ is invariant) and \eqref{eq-EscapeIn-W-inside-R1} (implying that, for all $i\in[N]$, $E_i(x)<0$ whenever $W(x)<0$ since $\nu_i\geq0$, so that $\Sub_s(W)\subseteq\nint\,\calR$ for all $s<0$).
Using this invariant region, we seek a quadratic function $V$ whose Lie derivative is smaller than $G(V)$ for all points outside the invariant region $\Sub_0(W)$.
This is enforced by \eqref{eq-EscapeIn-deriv} (implying that $\calL_AV(x)\leq G(V(x))$ whenever $W(x)\geq0$).
We also require that $\calR\subseteq\Sub_1(V)$.
This is enforced by \eqref{eq-EscapeIn-include-R0} (implying that, for all $i\in[M]$, $V(x)\leq1$ whenever $F_i(x)\leq0$ since $\lambda_i\geq0$).
Finally, we try to find the largest $r$ such that $\Sub_r(V)\subseteq\calR$.
This is enforced by \eqref{eq-EscapeIn-inside-R1} (implying that $\Sub_s(V)\subseteq\nint\,\calR$ for any $s<r$; similarly to the case of \eqref{eq-EscapeIn-W-inside-R1}).
See also Figure~\ref{fig-EscapeIn-simple} and Example~\ref{exa-EscapeIn-simple} for an illustration.

The significance of the above optimization program \eqref{eq-EscapeIn} is that any feasible solution provides an upper bound on the escaping time from $\calR$:

\begin{thm}\label{thm-EscapeIn-bound}
Assume $0\in\nint\,\calR$.
Let $(r,V,\ldots)$ be a feasible solution%
\footnote{Note that $r$ can be negative.}\!
of \eqref{eq-EscapeIn} with $G(V)=-1$.
Then, for any $x_0\in\calR$,
\[
t_*(x_0,\calR;A)\leq\max(V(x_0)-r,0)\leq\max(1-r,0).
\]
Similarly, let $(r,V,\ldots)$ be a feasible solution of \eqref{eq-EscapeIn} with $G(V)=-\gamma V$, such that $r>0$.
Then, for any $x_0\in\calR$,
\[
t_*(x_0,\calR;A)\leq\frac{\log^+\!\big(\frac{V(x_0)}r\big)}{2\gamma}\leq\frac{\log^+\!\big(\frac1r\big)}{2\gamma},
\]
where $\log^+(\alpha)\triangleq\log(\max(\alpha,1))$.
\end{thm}

\begin{pf}
Let $x_0\in\calR$ and denote $\tau_*=t_*(x_0,\calR;A)$.
Assume $\tau_*>0$, as otherwise the inequalities are trivially satisfied.
Let $x(\cdot)$ be the trajectory of the system $\dot{x}=Ax$, starting from $x_0$.
Then, for all $t\in[0,\tau_*]$, it holds that $W(x(t))\geq0$; indeed, otherwise $x(t)$ would be in an invariant set inside $\nint\,\calR$ (see ``Explanations'').
First, consider the model with $G(V)=-1$.
For all $t\in[0,\tau_*]$, it holds that $V(x(t))\geq r$ and $\smash{\dot{V}}(x(t))\leq-1$, by \eqref{eq-EscapeIn-include-R0}, \eqref{eq-EscapeIn-deriv} and the fact that $W(x(t))\geq0$.
This implies that
\[
r\leq V(x(\tau_*))=V(x_0)+\int_0^{\tau_*}\dot{V(x(t))}\,\diff t\leq V(x_0)-\tau_*.
\]
Hence, we find that $\tau_*\leq V(x_0)-r$.
Finally, the inequality $\tau_*\leq 1-r$ comes from the fact that $V(x_0)\leq1$ for every $x_0\in\calR$ (see ``Explanations'').

The proof for the model with $G(V)=-\gamma V$ is along the same lines: we get that, for all $t\in[0,\tau^*]$,
\[
V(x(t))\leq V(x_0)+\int_0^t-2\gamma V(x(s))\,\diff s.
\]
The classical argument (like Gr\"onwall's inequality) then implies that $r\leq V(x(\tau_*))\leq V(x_0)e^{-2\gamma\tau_*}$, concluding the proof.\hfill\qed
\end{pf}

It is not difficult to see that \eqref{eq-EscapeIn} with $G(V)=-1$ always admits a feasible solution.
Likewise, for any $\gamma<\sigma(A)$, \eqref{eq-EscapeIn} with $G(V)=-2\gamma V$ admits a feasible solution satisfying $r>0$.
In the numerical experiments (see Section~\ref{sec-comparisons}), we have used the value $\gamma=\sigma(A)/2$.
This choice comes from a trade-off between (i) flexibility in the shaping of the invariant regions (the constraint \eqref{eq-EscapeIn-deriv} is softer when $\gamma$ is smaller) and (ii) rate of decrease of $V$ along the trajectories of the system inside $\calR$ (equal to $e^{-2\gamma t}$).
From Theorem~\ref{thm-EscapeIn-bound}, it follows that by choosing $\gamma=\sigma(A)/2$, we have an upper bound that is in the worst case within a factor $2$ of the upper bound that we would get if using $\gamma$ close to $\sigma(A)$, while allowing much more flexibility in the shaping of the invariant regions, so that in practice the upper bound is better that the ones with $\gamma\approx\sigma(A)$.

Finally, the objective \eqref{eq-EscapeIn-obj} of \eqref{eq-EscapeIn} guarantees that the worst-case upper bound on $t_*(x_0,\calR;A)$, when $x_0$ varies in $\calR$, obtained from the model is the smallest possible within the considered framework (see Theorem~\ref{thm-EscapeIn-bound}).

\begin{exmp}\label{exa-EscapeIn-simple}
Consider the matrix $A = \left[\begin{smallmatrix} -1 & 3 \\ 0 & -1 \end{smallmatrix}\right]$, and the region $\calR$ depicted in Figure~\ref{fig-EscapeIn-simple}.
The application of \eqref{eq-EscapeIn} with $G(V)=-2\gamma V$ for this $A$ and this $\calR$ is illustrated in Figure~\ref{fig-EscapeIn-simple}.
We have used $\gamma=\sigma(A)/2=0.5$.
We have also represented a sample trajectory of the system, for $t\in[0,\tau]$ where $\tau=\log(V(x_0)/r)/(2\gamma)$ is obtained from the optimal solution of \eqref{eq-EscapeIn}; see Theorem~\ref{thm-EscapeIn-bound}.
We verify that the trajectory does not exit the region $\calR$ after time $\tau$.
\end{exmp}

\subsection{Upper bound on escaping time. Case II: $0\notin\calR$}\label{ssec-bound-escaping-time-outside}

Next, we consider the case where the origin is outside the region $\calR$.
Consider the following optimization program.

\emph{Decision variables:} a convex quadratic function $V$, and scalars $\lambda_1,\ldots,\lambda_M\geq0$, $\mu_1,\ldots,\mu_N\geq0$, $\nu_1,\ldots,\nu_N\geq0$ and $r\in\Re$;
\emph{Objective and constraints:}
\begin{subequations}\label{eq-EscapeOut}
\begin{talign}
&\!\!\makebox[10mm][l]{max}  r \label{eq-EscapeOut-obj}\\
&\!\!\makebox[10mm][l]{s.t.} V \leqq 1 + \lambda_i F_i, \quad \forall\,i\in[M], \label{eq-EscapeOut-include-R0}\\
&\!\!\makebox[10mm][l]{}     V \geqq r - \sum_{i=1}^N\mu_i E_i, \label{eq-EscapeOut-outside-R1}\\
&\!\!\makebox[10mm][l]{}     \calL_AV \leqq G(V) + \sum_{i=1}^N\nu_i E_i, \label{eq-EscapeOut-deriv}
\end{talign}
\end{subequations}
where $G(V)$ is equal to either $-1$ (leading to quadratically increasing bounds), or to $-2\gamma V$ for some fixed parameter $\gamma>0$ (leading to logarithmically increasing bounds).

\begin{figure}
\centering
\includegraphics[width=0.99\linewidth]{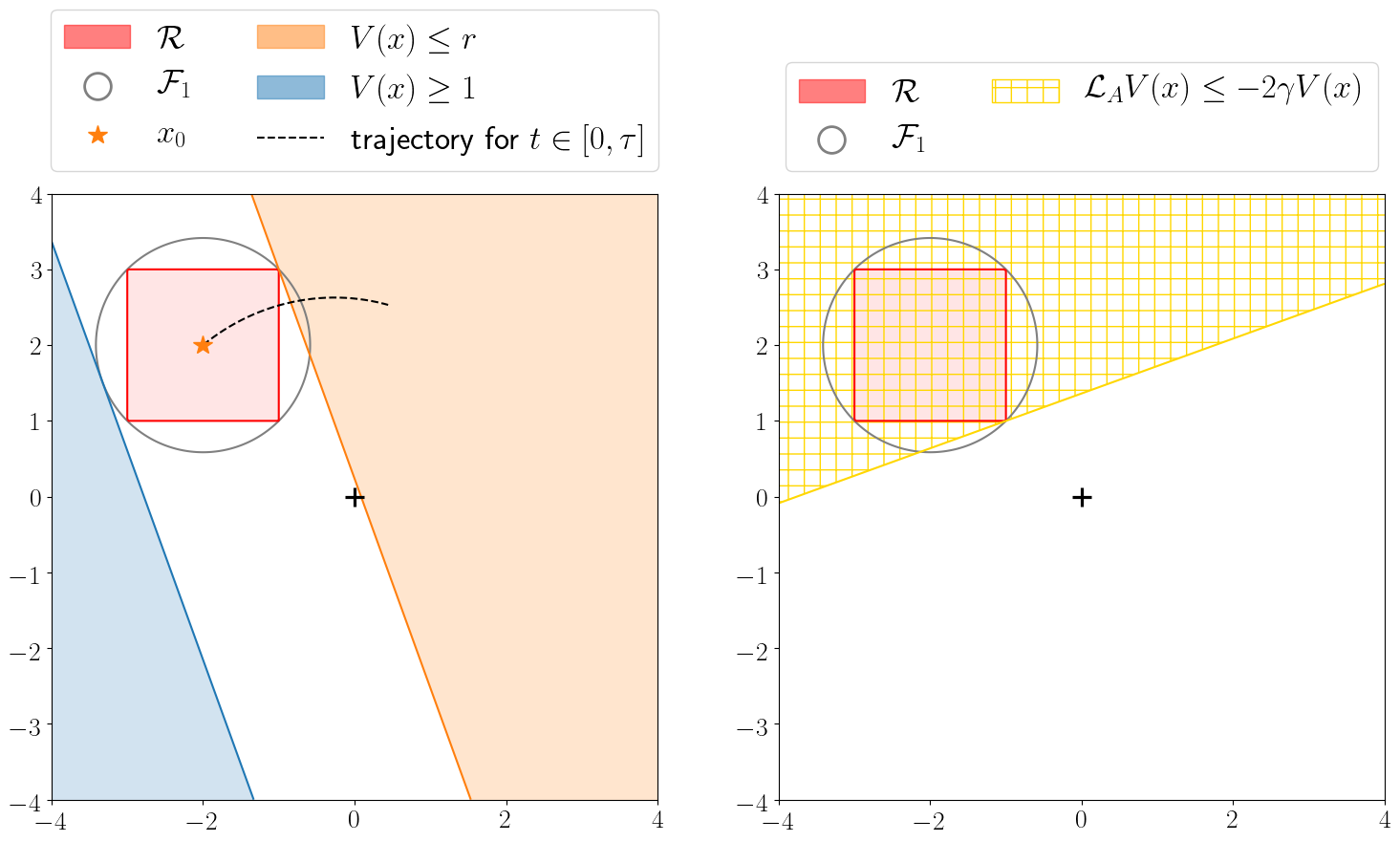}
\caption{Illustration of Example~\ref{exa-EscapeOut-simple}.
We verify that $\calR\subseteq\Sub_1(V)$, $\Sub_r(V)\cap\nint\,\calR=\emptyset$, and $\calL_A V(x)\leq G(V(x))$ for all $x\in\calR$.}
\label{fig-EscapeOut-simple}
\end{figure}

{\bfseries Explanations.} The goal is to find a quadratic function $V$ whose Lie derivative is smaller than $G(V)$ inside $\calR$.
This is enforced by \eqref{eq-EscapeOut-deriv} (implying that $\calL_AV(x)\leq G(V)$ whenever $E_i(x)\leq0$ for all $i\in[N]$ since $\nu_i\geq0$).
We also require that $\calR\subseteq\Sub_1(V)$.
This is enforced by \eqref{eq-EscapeOut-include-R0} (similar to Subsection~\ref{ssec-bound-escaping-time-inside}).
Finally, we try to find the largest $r$ such that $\Sub_r(V)$ is outside $\calR$.
This is enforced by \eqref{eq-EscapeOut-outside-R1} (implying that $V(x)\geq r$ whenever $E_i(x)\leq0$ for all $i\in[N]$ since $\mu_i\geq0$, so that $\Sub_s(V)\cap\calR=\emptyset$ for all $s<r$).
See also Figure~\ref{fig-EscapeOut-simple} and Example~\ref{exa-EscapeOut-simple} for an illustration.

Any feasible solution of the above optimization program provides an upper bound on the escaping time from $\calR$, in the exact same way as in Theorem~\ref{thm-EscapeIn-bound}:

\begin{thm}\label{thm-EscapeOut-bound}
Assume $0\notin\calR$.
Let $(r,V,\ldots)$ be a feasible solution%
\footnote{Note that $r$ can be negative.}\!
of \eqref{eq-EscapeOut} with $G(V)=-1$.
Then, for any $x_0\in\calR$,
\[
t_*(x_0,\calR;A)\leq\max(V(x_0)-r,0)\leq\max(1-r,0).
\]
Similarly, let $(r,V,\ldots)$ be a feasible solution of \eqref{eq-EscapeIn} with $G(V)=-\gamma V$, such that $r>0$.
Then, for any $x_0\in\calR$,
\[
t_*(x_0,\calR;A)\leq\frac{\log^+\!\big(\frac{V(x_0)}r\big)}{2\gamma}\leq\frac{\log^+\!\big(\frac1r\big)}{2\gamma},
\]
where $\log^+(\alpha)\triangleq\log(\max(\alpha,1))$.
\end{thm}

\begin{pf}
Similar to that for Theorem~\ref{thm-EscapeIn-bound}; omitted.\hfill\qed
\end{pf}

As for the model \eqref{eq-EscapeIn} in the previous subsection, it is not difficult to see that \eqref{eq-EscapeOut} with $G(V)=-1$ is always feasible.
The same holds with $G(V)=-2\gamma V$ and $r>0$, whenever $\gamma<\sigma(A)$.
In the numerical experiments (see Section~\ref{sec-comparisons}), we have used the value $\gamma=\sigma(A)/2$, which is a trade-off between (i) flexibility in the shaping of the invariant regions and (ii) rate of decrease of $V$ along trajectories of the system inside $\calR$ (see also Subsection~\ref{ssec-bound-escaping-time-inside} for details).

In the same way, the objective \eqref{eq-EscapeOut-obj} of \eqref{eq-EscapeOut} guarantees that the worst-case upper bound on $t_*(x_0,\calR;A)$, when $x_0$ varies in $\calR$, obtained from the model is the smallest possible within the considered framework (see Theorem~\ref{thm-EscapeOut-bound}).

\begin{exmp}\label{exa-EscapeOut-simple}
Consider the matrix $A = \left[\begin{smallmatrix} -0.1 & 1 \\ -1 & -0.1 \end{smallmatrix}\right]$, and the region $\calR$ depicted in Figure~\ref{fig-EscapeOut-simple}.
The application of \eqref{eq-EscapeOut}, with $G(V)=-2\gamma V$ for this $A$ and this $\calR$ is illustrated in Figure~\ref{fig-EscapeOut-simple}.
We have used $\gamma=\sigma(A)/2=0.05$.
We have also represented a sample trajectory of the system, for $t\in[0,\tau]$ where $\tau=\log(V(x_0)/r)/(2\gamma)$ is obtained from the optimal solution of \eqref{eq-EscapeOut}; see Theorem~\ref{thm-EscapeOut-bound}.
We verify that the trajectory exits the region $\calR$ before time $\tau$.
\end{exmp}

\section{Numerical experiments and comparisons}\label{sec-comparisons}

In this section, we would like to compare the performances of the two models of Section~\ref{sec-upper-bounds} (with $G(V)=-1$ and with $G(V)=-2\gamma V$).
We also compare them with the approach of \citet{rabi2020piecewise}, which relies on quadratic Lyapunov functions obtained by arbitrarily restricting the right-hand side of the Lyapunov equation, as in $A^\top\! P+PA=-I$.
To make the comparison, we have applied the different approaches on randomly generated matrices of dimension $n=10$, as explained below.

{\bfseries Random matrix generation.}
To generate random \emph{stable} matrices, we proceed as follows:
(1) Start from a random block-diagonal matrix $D$ with $2$-dimensional blocks on the diagonal.
Each block is obtained from a pair of real or complex conjugated eigenvalues $\lambda_1,\lambda_2$ defined as the roots of $\lambda^2+\sqrt{\alpha}\lambda+\beta/4$ where $\alpha$ and $\beta$ are chosen uniformly at random over the interval $[0,1]$.%
\footnote{This implies that $\lambda_1,\lambda_2$ have negative real parts with probability $1$.
Moreover, they are real with probability $\frac12$ (useful to investigate the effect of both real and complex eigenvalues on the performance of the techniques).}\!
(2) From $D$, build the random matrix $A$ by applying the similarity transformation $A=UDU^{-1}$, where the elements of $U$ are randomly generated from the standard normal distribution (which implies that $U$ is invertible with probability $1$).

{\bfseries Regions definition.}
We also have to define the region $\calR$.
We made the choice to consider a fixed region because possible transformations are included in the random similarity transformations $UDU^{-1}$.
For the case I (origin inside), we let $\calR$ be the set of points $x$ such that $-2.5\leq x[1]\leq1.5$ and $-2\leq x[i]\leq2$ for $i=2,\ldots,10$ (where $x[i]$ is the $i$th component of $x$).
We also let $x_0=[-1.5,-1,\ldots,-1]^\top$.
For the case II (origin outside), we let $\calR$ be the set of points $x$ such that $-6\leq x[i]\leq-2$ for $i=1,2$, and $-2\leq x[i]\leq2$ for $i=3,\ldots,10$.
We let $x_0=[-5,-5,-1,\ldots,-1]^\top$.

We have used the SDP solver Mosek (https://www.mosek .com) with the package JuMP \citep{dunning2017jump} in Julia, to compute the solutions of the different models.
All computations were done on a laptop with processor Intel Core i7-7600u and 16 GB RAM running Windows.

\subsection{Experimental results}

\begin{figure}
\centering
\includegraphics[width=0.9\linewidth]{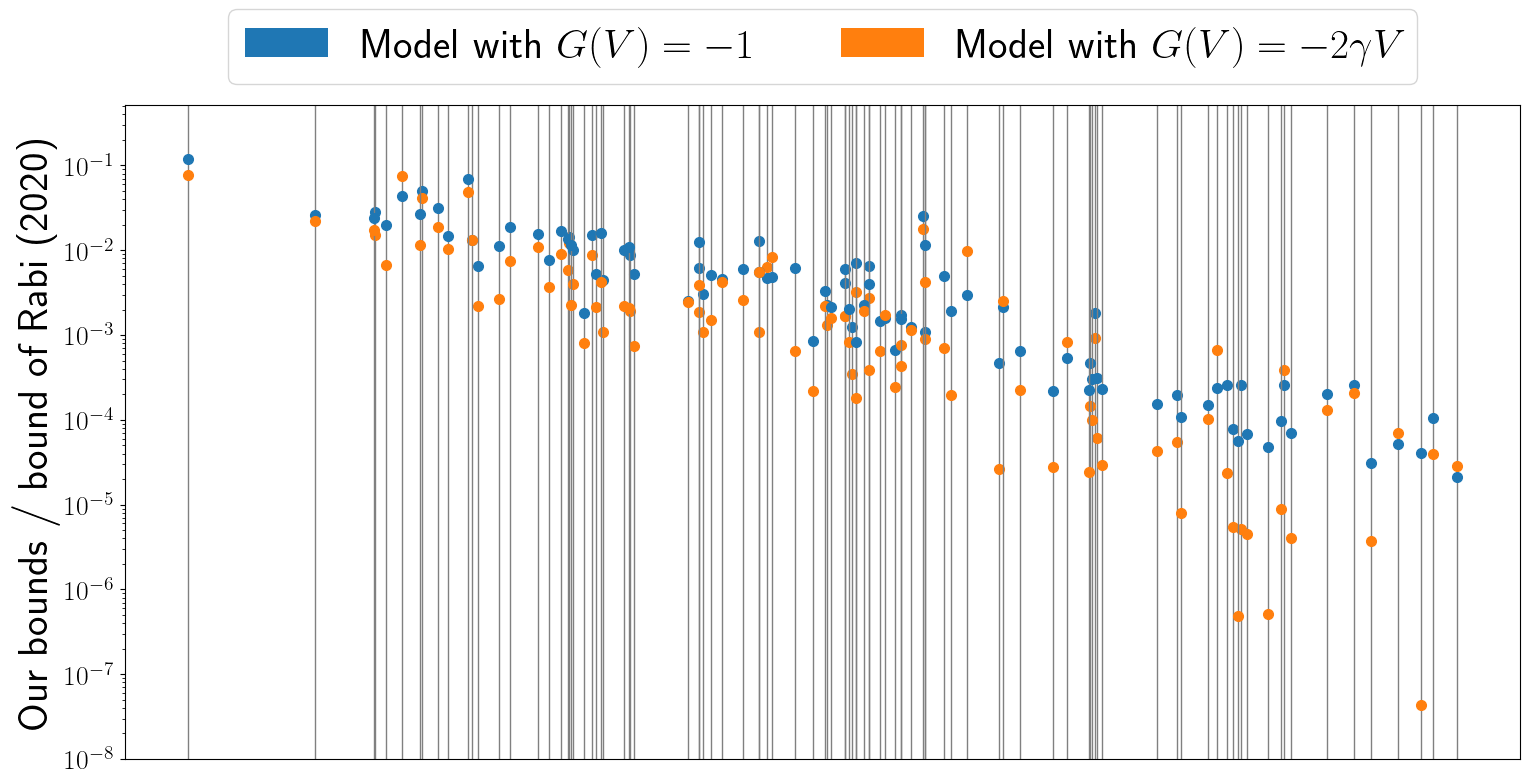}\\
\includegraphics[width=0.9\linewidth]{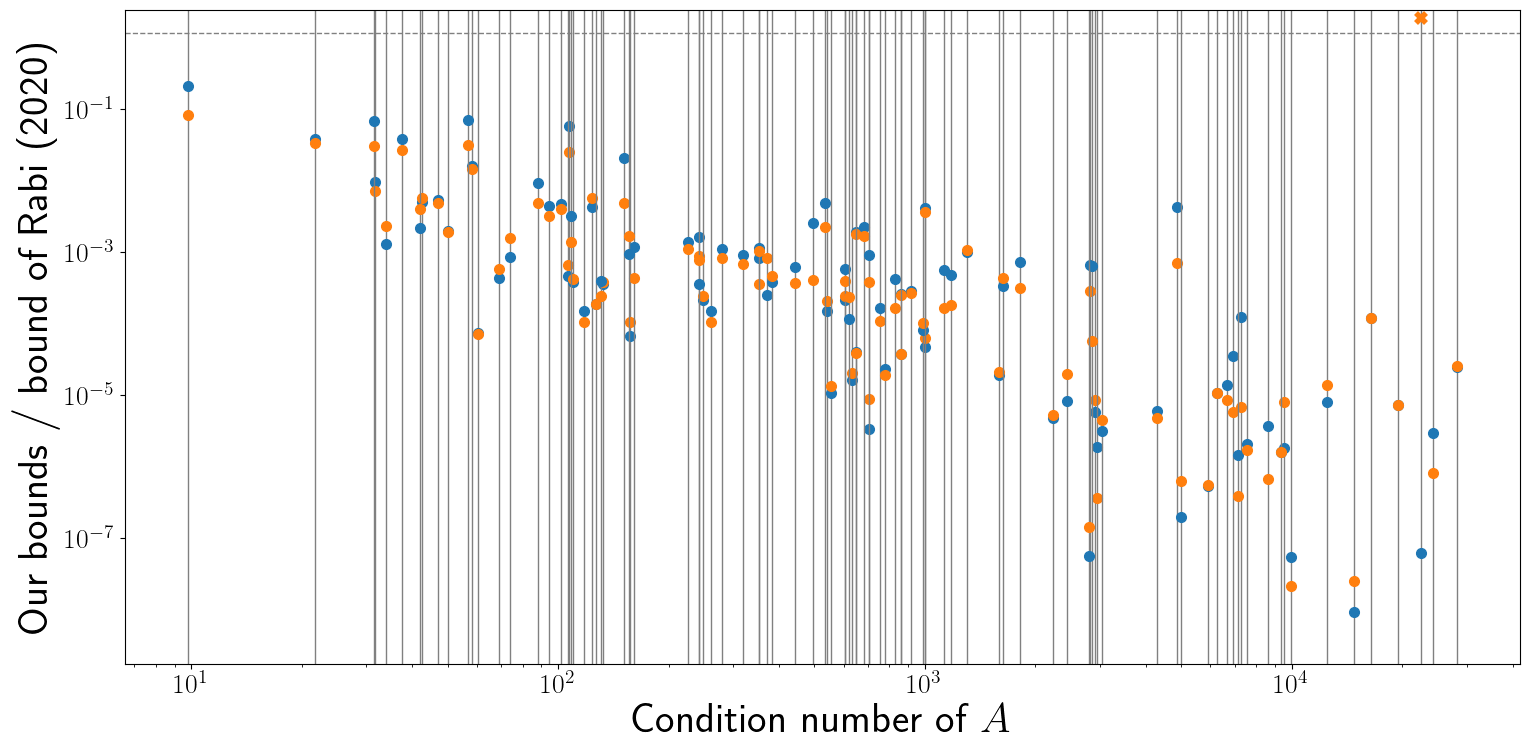}
\caption{\emph{Top:} Case I (origin inside).
\emph{Bottom:} Case II (origin outside).
The same matrices were used for both cases.}
\label{fig-CompareBounds}
\end{figure}

\begin{figure}
\centering
\includegraphics[width=0.9\linewidth]{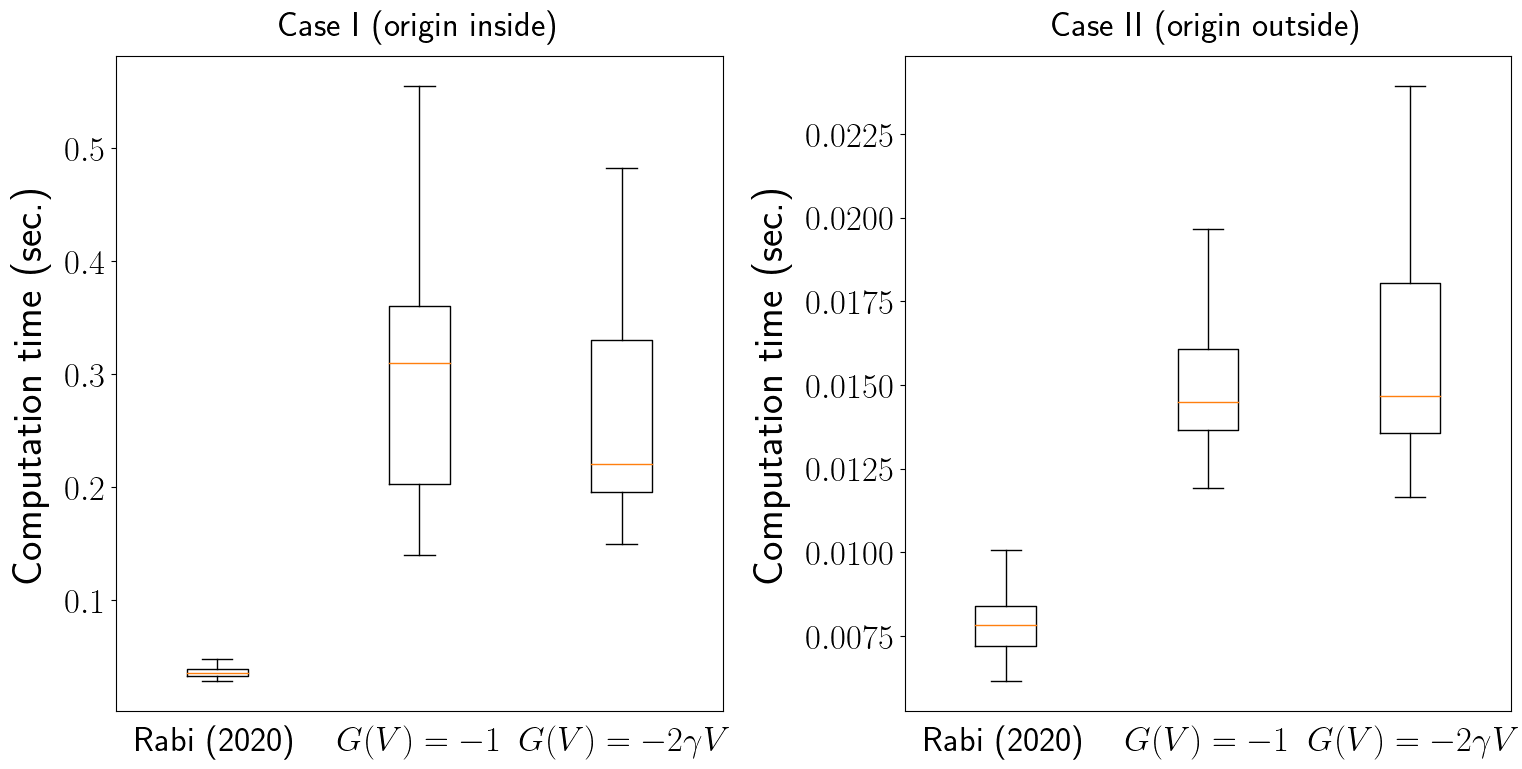}
\caption{Box plots of computation times of the different approaches.}
\label{fig-computation-times}
\end{figure}

For $100$ randomly generated matrices $A_i$, $i=1,\ldots,100$, we have used the approaches described in Subsections~\ref{ssec-bound-escaping-time-inside} and~\ref{ssec-bound-escaping-time-outside} and in \citet{rabi2020piecewise} to compute upper bounds on the escaping time $t_*(x_0,\calR;A_i)$.
The plots in Figure~\ref{fig-CompareBounds} show the pairs $(\kappa_i,\rho_i)$, where $\kappa_i$ is the condition number of $A_i$, and $\rho_i$ is the ratio between the upper bound obtained from the models presented in this paper (with different colors for the different models) divided by the upper bound obtained using the approach of \citet{rabi2020piecewise}.%
\short{\footnote{\short{Comparisons between the obtained upper bound and the actual escaping time $t_*(x_0,\calR;A_i)$ for the different models presented in this paper are discussed in the extended version of this work \citep[see][]{berger2021bounds}.}}\!\textsuperscript{,}\!}%
\footnote{In one case, the solver failed to compute a feasible solution within the allotted time (using the solver's default options).
This case is shown in the second chart of Figure~\ref{fig-CompareBounds}; it is indicated by the cross near the chart's top-right corner.}
The goal is to evaluate the improvement provided by the approaches in this paper, and how such improvements depend on the condition number of $A$.
Finally, the computation times of the different approaches are given in Figure~\ref{fig-computation-times}.
\extended{For the interested reader, the ratio between the obtained upper bound and the actual escaping time $t_*(x_0,\calR;A_i)$ for the different models of this paper are given Figure \ref{fig-CompareBounds-absolute}.}

\extended{%
\begin{figure}
\centering
\includegraphics[width=0.9\linewidth]{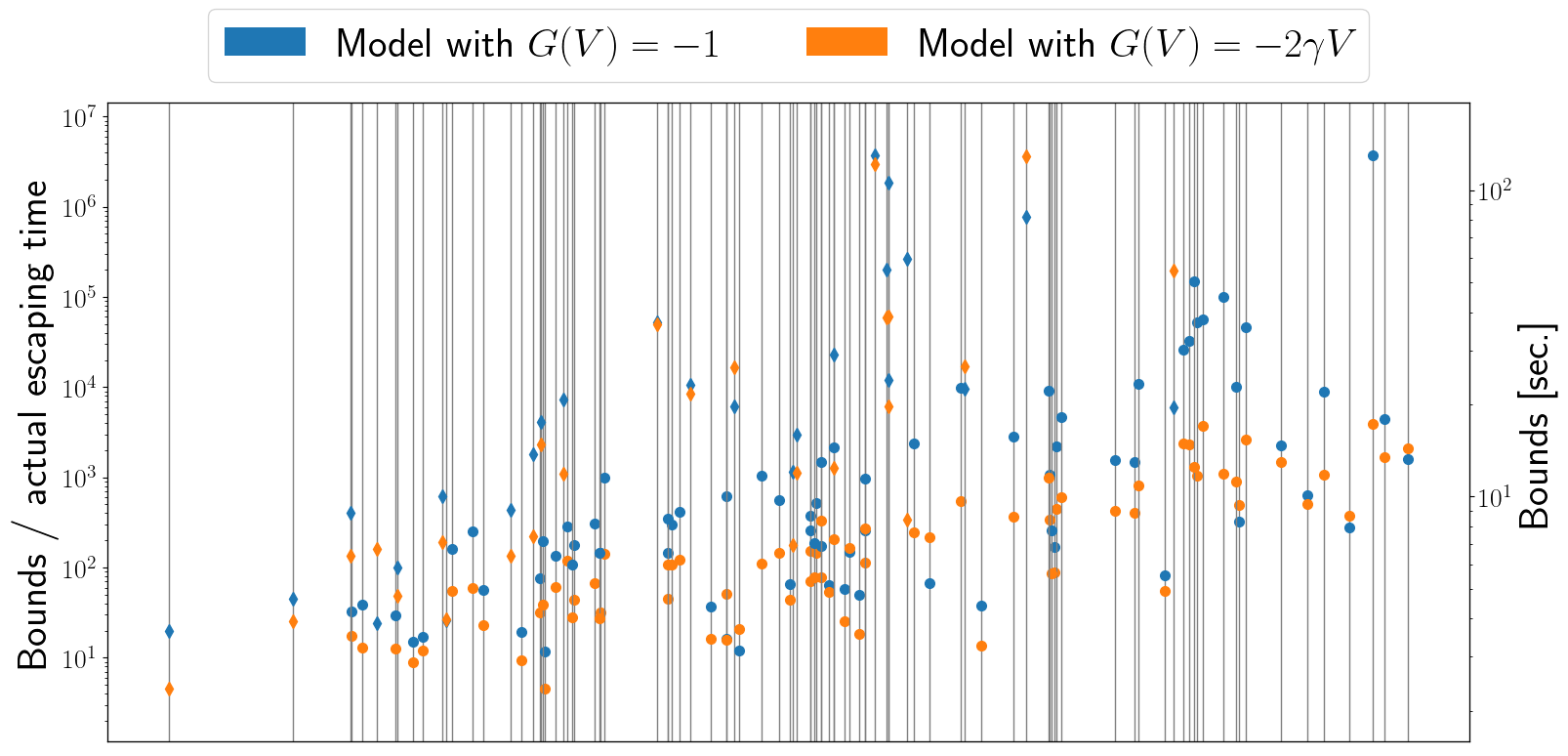}\\
\includegraphics[width=0.9\linewidth]{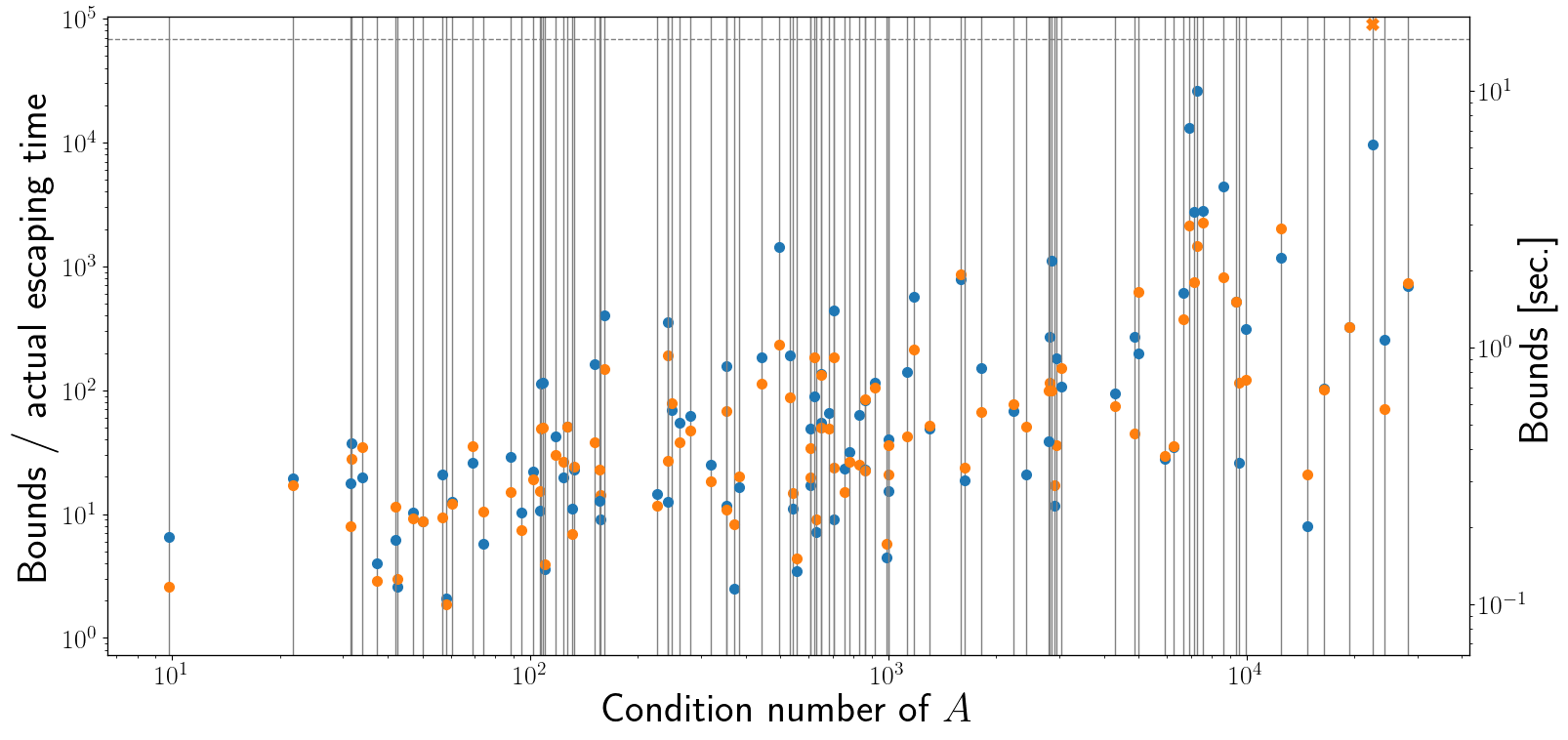}
\caption{\extended{\emph{Top chart:} Case I (origin inside).
Since $0\in\nint\,\calR$, it may happen that $x(\cdot)$ never exits $\calR$, so that $t_*(x_0,\calR;A_i)=0$.
In this case, two situations can occur: either the upper bound provided by the optimization model is zero, or this bound is positive.
In the first situation, the ratio between the upper bound (equal to zero) and the actual escaping time (equal to zero) is taken as one, by our convention.
In the second situation, the ratio is infinite, and therefore, we have represented the point $(\kappa(A_i),\tau_i)$ instead of $(\kappa(A_i),\rho_i)$, where $\tau_i$ is the upper bound provided by the optimization model.
Such points are marked as diamonds, to be distinguishable from the usual points (represented by circles).  For any diamond-marked point, the value of $\tau_i$ is indicated on the vertical axis on the right.
\emph{Bottom chart:} Case II (origin outside).
The same set of matrices $A$ were used for both cases.}}
\label{fig-CompareBounds-absolute}
\end{figure}}

\subsection{Conclusions from the numerical experiments}

In Figure~\ref{fig-CompareBounds}, we observe that the models presented in this paper provide upper bounds on the escaping time that are between $10^1$ and $10^8$ times better than the upper bounds obtained with the approach described in \citet{rabi2020piecewise}.
We also notice that the higher the condition number of $A$, the better the improvement tends to be.
Finally, Figure~\ref{fig-CompareBounds} also shows that the models with $G(V)=-2\gamma V$ provide slightly better upper bounds than the models with $G(V)=-1$ for the case I, and comparable upper bounds for the case II.

Regarding computation time, the approach of \citet{rabi2020piecewise} is in average up to $10$ times faster than the methods in this paper.
This is not really surprising since our models involves a number of variables growing as $n^2+N+M$, while the approach of \citet{rabi2020piecewise} involves a number of variables proportional to $N$ (number of ellipsoids $\calE_i$).

\section{Conclusions}

In this paper, we improved an existing way of using quadratic Lyapunov theory to compute upper bounds on the set escaping time of stable linear systems.
Our bounds are based on convex programming techniques.
In particular, our approach uses Lyapunov-like functions that are not necessarily centered at the origin, and whose rates of decrease are guaranteed only in a specific subset of the state space.
Numerical experiments demonstrate that in many cases our approach is highly favorable in terms of the accuracy of the upper bound (up to $10^8$ times better), while requiring only a modest increase in computation time (up to $10$ times slower).

For further work, we plan to investigate the usefulness of both approaches in practical applications, for instance as part of a stand-alone software for the simulation of switched linear systems (see the introduction).
Finally, we also plan to investigate ways to speed up the computation of the Lyapunov-like functions in our models; for instance, by leveraging the fact that in general we do not need an accurate optimal solution of the optimization models (a near-optimal feasible solution is often sufficient), or that a feasible initial solution can always be obtained from the Lyapunov equation.

\bibliography{myrefs}






\end{document}